\documentclass[12pt,twoside,final]{amsart}
\usepackage[psamsfonts]{amssymb}
\usepackage{times,a4wide}
\usepackage{xcolor}

\usepackage{ulem}
\usepackage{amsrefs}
\usepackage{hyperref}

\theoremstyle{plain}
\newtheorem*{theorem*}{Theorem}

\newtheorem{theorem}{Theorem}[section]
\newtheorem{proposition}[theorem]{Proposition}
\newtheorem*{maintheorem*}{Main Theorem}
\newtheorem*{proposition*}{Proposition}
\newtheorem{corollary}[theorem]{Corollary}
\newtheorem*{corollary*}{Corollary}
\newtheorem{lemma}[theorem]{Lemma}
\newtheorem*{lemma*}{Lemma}

\theoremstyle{definition}
\newtheorem{remark}[theorem]{Remark}

\newtheorem*{remark*}{Remark}
\newtheorem*{remarks*}{Remarks}
\newtheorem*{conjecture*}{Conjecture}
\theoremstyle{definition}

\newenvironment{poliabstract}[1]
  {\renewcommand{\abstractname}{#1}\begin{abstract}}
  {\end{abstract}}

\allowdisplaybreaks

\usepackage[active]{srcltx}
\usepackage{tikz}
\usepackage{pgfplots}
\usepackage{float}
\RequirePackage{pgfplots}
\usetikzlibrary{shadows}
\pgfplotsset{compat=1.13} 

\newcommand{\C}{\mathbb{C}}

\newcommand{\N}{\mathbb{N}}

\newcommand{\E}{\mathbb{E}}

\renewcommand{\P}{\mathbb{P}}

\newcommand{\Var}{\operatorname{Var}}

\title[Separated determinantal point processes] 
{Separated determinantal point processes and generalized Fock spaces} 

\author[G. Lamberti]{Giuseppe Lamberti}

\address{Univ. Bordeaux, CNRS, Bordeaux INP, IMB, UMR 5251, F-33400, Talence, France}

\author[X. Massaneda]{Xavier Massaneda}

\address{Departament de Matem\`atiques i Inform\`atica,
Universitat  de Barcelona, Gran Via 585, 08007-Bar\-ce\-lo\-na, Catalonia}

\thanks{The first author partially supported by the project UBGRS 2.0 (ANR-20-SFRI-0001). Second author partially supported by the Generalitat de Catalunya (grant 2021 SGR 00087) and the spanish Ministerio de Ciencia e Innovaci\'on (project PID2021-123405NB-I00).}

\date{\today}


\begin{document}

\begin{poliabstract}{Abstract} 
Let $\Lambda_\phi$ denote the determinantal point processes associated to the generalized Fock space defined by a doubling subharmonic weight $\phi$. We provide a characterization of the processes $\Lambda_\phi$ which almost surely form a separated sequence in $\C$. Furthermore, we emphasize the role of intrinsic repulsion in determinantal processes by comparing $\Lambda_\phi$ with the Poisson process of the same first intensity. As an application, we show that
the determinantal process $\Lambda_\alpha$ associated to the canonical weight $\phi_\alpha(z)=|z|^\alpha$, $\alpha>0$, is almost surely separated if and only if $\alpha<4/3$. In contrast, the Poisson process $\Lambda_\alpha^P$ having the same first intensity as $\Lambda_\alpha$ is almost surely separated if and only if $\alpha<1$.
\end{poliabstract}

\maketitle

\section{Introduction}\label{sec:intro}

In this paper we provide conditions so that a particular family of determinantal point processes in the complex plane are almost surely separated.  Recall that a sequence $\Lambda=\{\lambda_k\}_{k\geq 1}\subset\C$ is \textit{separated} if
\[
 \inf_{j\neq k} |\lambda_k-\lambda_j|>0.
\]

Separation plays an important role in many numerical and function theoretic problems, such as the description of interpolating and sampling sequences for various spaces.

Intuitively, a (simple) point process is a random configuration of points, but it is more convenient to think of it as a random measure of the form
\[
 \tau_\Lambda=\sum_{\lambda\in \Lambda}\delta_\lambda\ ,
\]
where $\Lambda$ is a finite or countable subset of $\C$. The distribution (or law) of the point process is then determined by the random variables
\[
N(B)=\#(\Lambda\cap B)=\int_B d\tau_\Lambda,\qquad \text{$B\subset\C$ compact}.
\]
For a background on random measures and point processes we refer the reader to \cite{Da-Ve}.

We start by recalling the definition of the generalized Fock space $\mathcal F_{\phi}$ defined by a doubling weight $\phi$.
Let $\phi$ be a subharmonic function in $\C$ with doubling Laplacian $\nu=\Delta\phi$, i.e., for which there exists $C>0$ such that 
\[
 \nu(D(z,2r))\leq C\, \nu(D(z,r)) ,\qquad z\in\C,\quad r>0.
\]
Then, denote by $\mathcal F_{\phi}$ the Fock space of entire functions $f$ such that
\[
\int_{\C}|f(z)|^2 e^{-2\phi(z)}\Delta\phi(z)<+\infty.
\]
This is a reproducing kernel Hilbert space whose kernel $K_\phi$ is also referred to as the Bergman kernel of $\mathcal F_{\phi}$ (see next section for precise definitions).

The focus of this paper is on the determinantal point processes $\Lambda_\phi$ associated to $K_\phi$. We proceed to explain this in further detail. According to a theorem of Macch\`i and Soshnikov (\cite{soshnikov00}*{Theorem 3}, see also \cite{HKPV}*{Lemma 4.5.1}) whenever an Hermitian kernel $K(z,\zeta)$ defines a self-adjoint operator $\mathcal K$ on $L^2(\C,\mu)$ which is locally trace class with all eigenvalues in $[0,1]$, there exists a determinantal point process $\Lambda$ associated with $K$ and $\mu$. This means that for any collections of disjoint sets $B_1,\dots, B_n\subset\C$ one has
\[
 \mathbb E\left[\prod_{k=1}^n \#(B_k\cap\Lambda)\right]=
 \int_{B_1}\cdots \int_{B_n}\det\bigl(K(z_i,z_j)\bigr)_{1\leq i,j\leq n} d\mu(z_1)\cdots d\mu(z_n).
\]
The integrand is usually called the $n^{\text{th}}$-\textit{correlation function} of $\Lambda$. For the particular case corresponding to $n=1$, this is
\[
 \mathbb E[N(B)]=\int_{B} K(z,z)\, d\mu(z),
\]
and the measure $K(z,z)\, d\mu(z)$ is called the \textit{first intensity} (or average distribution) of the process. 

The Bergman kernel $K_\phi$ defines a projection operator from $L^2(\C,\mu_\phi)$ onto $\mathcal F_\phi$, so it satisfies the hypothesis of the aforementioned theorem by Macch\`i and Soshnikov (the definition of $\mu_\phi$ is provided in the following section). Consequently, the determinantal point process $\Lambda_\phi$ associated to $K_\phi$ and $\mu_\phi$ is well-defined.

The processes we study can be seen as generalizations of the process studied in \cite{Bufetov17}, where the authors investigate the Palm measures associated to $\mathcal F_\phi$, under the assumption that $\Delta\phi$ is uniformly bounded both above and below.

We mention that determinantal processes were first considered by Macch\`i \cite{Ma} as a model for fermi\-ons, and since then have appeared naturally in a surprising variety of problems (see e.g. \cite{HKPV}*{Chapter 4}).

Our main result, Theorem~\ref{main}, characterizes the processes $\Lambda_\phi$ which are almost surely separated. As a by-product of its proof, we notice that when $\Lambda_\phi$ is a.s. separated, then its upper density is $0$. This automatically provides a characterization of the processes $\Lambda_\phi$ which are almost surely interpolating for the classical Fock spaces $\mathcal F_{\beta|z|^2}$, $\beta>0$ (see Corollary~\ref{cor:interpolation}). Of course, the analogous result holds as well for the $L^p$ versions of these classical Fock spaces, or for the Hilbert spaces of entire functions investigated in \cites{BOC,OCS}.

A distinctive feature of determinantal processes is their built-in repulsion (i.e. a ne\-gative correlation of nearby points), which can be quantified by the decay of the kernel. To assess the role played by this repulsion, we compare $\Lambda_\phi$ with the Poisson process $\Lambda_\phi^P$ of ``totally independent'' points having the same first intensity (average distribution) as $\Lambda_\phi$.
The existence of such a Poisson process is established, for instance, in \cite{LaPe}*{Section 3.2}.
An adaptation of \cite{HaMa}*{Theorem 1.2} shows that, by the big fluctuations of the Poisson process, 
much stricter conditions are required for $\Lambda_\phi^P$ to be almost surely separated (Theorem~\ref{thm:Poisson}).

The plan of the paper is as follows. In
Section~\ref{sec-det} we gather the main properties of the Fock spaces $\mathcal F_\phi$ and their associated processes $\Lambda_\phi$. We also state the main results. The proof of the main result (Theorem~\ref{main}) and its consequences is given in Section~\ref{sec:main}. The key step in this proof is the estimate of the probability that two or more points fall in a cell of fixed size (Proposition~\ref{prop:traces}).
Finally, Section~\ref{sec:Poisson} is devoted to briefly recall the proof of Theorem~\ref{thm:Poisson}. 

A final word about notation: the expression $A\preceq B$ means that there
exists a constant $C>0$, independent of whatever arguments are involved, such
that $A\leq C B$. If both $A\preceq B$ and $B\preceq A$ then we write $A\simeq
B$.

{\bf Acknowledgements:} We would like to thank Jordi Marzo and Joaquim Ortega-Cerd\`a for their help on determinantal processes and for many valuable discussions.

\section{Determinantal processes induced by generalized Fock spaces. Results.}\label{sec-det}
Let $\phi$ be a subharmonic function with doubling Laplacian and let $\mathcal F_\phi$ be the associated Fock space. In order to express the norm in $\mathcal F_\phi$ it is more convenient to work with a suitable regularization of $\nu=\Delta\phi$.
Following the ideas of M. Christ \cite{Christ}, for each $z\in\C$ denote by $\rho(z)$ (or $\rho_\phi(z)$ if we want to stress the dependency on $\phi$) the positive radius such that
\[
 \nu(D(z,\rho(z))=1.
\]
The function $\rho^{-2}$ can then be seen as a regularized version of $\Delta\phi$, as described in \cite{Christ}.

It is not difficult to see (e.g. \cite{MMOC}*{p.869}) that
\begin{equation}\label{eq:rho-lip}
 |\rho(z)-\rho(\zeta)|\leq |z-\zeta|\, ,\qquad z,\zeta\in \C.
\end{equation}
In particular, $\rho(z)$ is continuous.

Define
\[
 \mathcal F_\phi=\Bigl\{ f\in H(\C) : \|f\|_\phi^2:=c_\phi\int_{\C} |f(z)|^2 e^{-2\phi_(z)} \, \frac{dm(z)}{\rho^2(z)}<+\infty\Bigr\}.
\]
The constant $c_\phi$ is chosen so that $\|1\|_{\phi}=1$; hence
\begin{equation}\label{eq:normalization_mu}
 d\mu_{\phi}(z):=c_\phi\, e^{-2\phi(z)} \frac{dm(z)}{\rho^2(z)}
\end{equation}
is a probability measure in $\C$.

Equipped with the scalar product
\[
 \langle f,g \rangle_{\phi}=\int_{\C} f(z)\, \overline{g(z)}\, d\mu_{\phi}(z),\qquad f,g\in \mathcal F_{\phi},
\]
the space $\mathcal F_{\phi}$ is a reproducing kernel Hilbert space whose kernel is denoted by $K_\phi$. Specifically, this means that $K_\phi$ is holomorphic in $z$, anti-holomorphic in $\zeta$, and for any $f\in \mathcal F_\phi$,
\[
 f(z)=\langle f,\overline{K(z,\cdot)}\rangle_{\phi}.
\]

Let now $\Lambda_\phi$ be the determinantal point process associated to the Hermitian kernel $K_\phi(z,\zeta)$ and the reference measure $\mu_{\phi}$. Its first intensity is, by definition,
\[
 K_\phi (z,z)\, d\mu_\phi(z).
\]
According to Lemma 21 in \cite{MMOC}; there exists $C>0$ such that
\begin{equation}\label{eq:kernel-diagonal}
C^{-1} e^{2\phi(z)}\leq K_\phi(z,z)\leq C\, e^{2\phi(z)}\  \qquad z\in\C.
\end{equation}
Therefore, by \eqref{eq:normalization_mu}, the first intensity of $\Lambda_\phi$ is comparable to the regularization of the Laplacian of $\phi$: 
\begin{equation}\label{eq:1st-intensity}
 K_\phi (z,z)\, d\mu_\phi(z)\simeq \frac{dm(z)}{\rho^2(z)}.
\end{equation}
In particular, given any Borel set $B\subset \C$, the expected value of the counting variable $N_\phi(B)=\#(\Lambda_\phi\cap B)$ is
\begin{equation}\label{eq:EB}
 \mathbb E[N_\phi(B)]=\int_B K_\phi (z,z)\, d\mu_\phi(z)\simeq \int_B\frac{dm(z)}{\rho^2(z)}.
\end{equation}

We are ready to state our main result. 

\begin{theorem}\label{main} Let $\phi$ be a doubling subharmonic function and let $\Lambda_\phi$ be its associated determinantal point process. Then
\[
 \P\bigl(\Lambda_\phi\ \text{is separated}\bigr)=
 \begin{cases}
  1\quad \text{if $\displaystyle{\int_{\C}}\dfrac{dm(z)}{\rho^6(z)}<+\infty$},\\
  0\quad\text{if $\displaystyle{\int_{\C}}\dfrac{dm(z)}{\rho^6(z)}=+\infty$}.
 \end{cases}
\]
\end{theorem}

\begin{remark}\label{rem:alternative}
According to \cite{MMOC}*{Theorem 14}, there always exist $\psi\in \mathcal C^\infty(\C)$ subharmonic and $C>0$ such that $\Delta\psi$ is a doubling measure, $|\psi-\phi|\leq C$ and $\Delta\psi\simeq \rho_{\psi}^{-2}\simeq  \rho_{\phi}^{-2}$, 
so that $\mathcal F_\phi=\mathcal F_\psi$ and $\|f\|_\phi\simeq \|f\|_\psi$ for all $f\in\mathcal F_\phi$.
In particular, the first intensity of $\Lambda_\phi$ and $\Lambda_\psi$ are both comparable to $\Delta\psi\, dm$.
The proof of Theorem~\ref{main} will show that $\Lambda_ \phi$ is separated if and only if $\Lambda_\psi$ is separated, 
thus we could as well assume from the beginning that $\phi$ is regular and rewrite
the integral above as
\[
 \int_{\C} (\Delta\phi(z))^3 dm(z).
\]
However, we prefer to keep the form given in the statement, because that is how it appears in the proof.
\end{remark}

Canonical examples of the weights considered in Theorem~\ref{main} are $\phi_\alpha(z)=|z|^\alpha$, $\alpha>0$. For $\alpha=2$, the distribution of $\Lambda_\alpha:=\Lambda_{\phi_\alpha}$ is invariant by translations (see e.g. \cite{RiVi}). In particular, for any $\varepsilon>0$ the probabilities $\P[\#(\Lambda_2\cap D(z,\varepsilon))]$ are strictly positive and independent of $z\in \C$. Then the arguments of Bufetov, Qiu and Shamov for the hyperbolic disk \cite{BQS}*{Lemma 1.3} apply similarly to $\C$ to conclude that $\Lambda_\alpha$ is a.s. not se\-parated. 

Since $\Delta|z|^\alpha=\alpha^2 |z|^{\alpha-2}$, Theorem~\ref{main} in this particular case yields the following.

\begin{corollary}\label{cor:alpha} Let $\Lambda_\alpha$ be the determinantal point process associated to the weight $\phi_\alpha(z)=|z|^\alpha$, $\alpha>0$. Then, 
\[
 \P\bigl(\Lambda_\alpha\ \text{is separated}\bigr)=
 \begin{cases}
  1\quad \text{if $\alpha<4/3$}\\
  0\quad\text{if $\alpha\geq 4/3$}.
 \end{cases}
\]
\end{corollary}

The proof of Theorem~\ref{main} shows that, in case $\Lambda_\phi$ is separated, its upper density
\[
 D_+(\Lambda_\phi):=\limsup_{r\to\infty}\sup_{z\in\C}\frac{\#(\Lambda_\phi\cap D(z,r))}{\pi r^2}
\]
is zero. By a well-known theorem by K. Seip and R. Wallstén \cites{Seip92a,Seip92b}, a sequence $\Lambda$ in the plane is interpolating for the classical Fock space $\mathcal F_{\beta|z|^2}$, $\beta>0$, if and only if it is separated and $D_+(\Lambda)<\beta$ (see e.g \cite{Zhu}*{Chapter 3}).

\begin{corollary}\label{cor:interpolation} Let $\phi$ be a doubling subharmonic weight and let $\Lambda_\phi$ be its associated determinantal process. Then, for any $\beta>0$,
\[
 \P\bigl(\Lambda_\phi\ \text{is $\mathcal F_{\beta|z|^2}$-interpolating}\bigr)=
 \begin{cases}
   1\quad \text{if $\displaystyle{\int_{\C}}\dfrac{dm(z)}{\rho^6(z)}<+\infty$}\\
  0\quad\text{if $\displaystyle{\int_{\C}}\dfrac{dm(z)}{\rho^6(z)}=+\infty$}.
 \end{cases}
\]
\end{corollary}
The results in \cite{BOC} and \cite{OCS} characterize interpolating sequences for weighted Fock spaces whith Laplacian bounded both above and below. The characterization is given in terms of two conditions: separation and an appropriate density bound. Consequently, Corollary~\ref{cor:interpolation} remains valid if $\mathcal F_{\beta|z|^2}$ is replaced with the Fock spaces studied in these references.

We finish this section by comparing the results above with the analogues for the Poisson process with the same first intensity. Consider 
\[
 d\sigma_\phi(z)=\frac{dm(z)}{\rho^2(z)}
\]
and denote by $\Lambda_\phi^P$ the Poisson process with underlying measure $\sigma_\phi$. Observe that, by \eqref{eq:1st-intensity},  
$\Lambda_{\phi}^P$ can be seen as the Poisson point process having the same average distribution as the determinantal process $\Lambda_{\phi}$.

The Poisson process $\Lambda_\phi^P$ is characterized by two properties:
\begin{itemize}
 \item [(i)] for any Borel set $A\subset\C$ the counting function $N(A)=\#(A\cap \Lambda_{\phi}^P)$ has Poisson distribution of parameter $\sigma_{\phi}(A)$,
 \item [(ii)] for two disjoint sets $A,B\subset\C$ the variables $N(A)$, $N(B)$ are independent. 
\end{itemize}

A straightforward adaptation of Theorem 1.2 in \cite{HaMa} provides the following characterization.

\begin{theorem}\label{thm:Poisson} 
Let $\phi$ be a doubling subharmonic weight and let $\Lambda_\phi^P$ be the Poisson process with underlying measure $\sigma_\phi$. Then 
\[
 \P\bigl(\Lambda_\phi^P\ \text{is separated}\bigl)=
 \begin{cases}
  1\quad \text{if $\displaystyle{\int_{\C}}\dfrac{dm(z)}{\rho^4(z)}<+\infty$}\\
  0\quad\text{if $\displaystyle{\int_{\C}}\dfrac{dm(z)}{\rho^4(z)}=+\infty$}.
 \end{cases}
\]
\end{theorem}

For the cases $\phi(z)=|z|^\alpha$, this means that $\Lambda_\phi^P$ is almost surely separated if and only if $\alpha<1$.

As in the determinantal case, we have an immediate consequence for classical Fock interpolation (or for the more general spaces studied in \cites{BOC,OCS}).

\begin{corollary}\label{cor:interpolation-P} Let $\phi$ be a doubling subharmonic weight and let $\Lambda_\phi^P$ be the Poisson process with underlying measure $\sigma_\phi$. Then, for any $\beta>0$,
\[
 \P\bigl(\Lambda_\phi\ \text{is $\mathcal F_{\beta|z|^2}$-interpolating}\bigr)=
 \begin{cases}
   1\quad \text{if $\displaystyle{\int_{\C}}\dfrac{dm(z)}{\rho^4(z)}<+\infty$}\\
  0\quad\text{if $\displaystyle{\int_{\C}}\dfrac{dm(z)}{\rho^4(z)}=+\infty$}.
 \end{cases}
\]
\end{corollary}

\section{Determinantal process: proof of Theorem~\ref{main} and Corollary \ref{cor:interpolation}}\label{sec:main}

\subsection{Preliminaries}\label{subsec:prelim} 
We start this section with some technical, well-known properties of $\mathcal F_\phi$.

A first elementary property to be used throughout the paper is the following (see \cite{MMOC}*{Corollary 3}): for any $r>0$ there exists $c_r>0$ such that
\begin{equation}\label{eq:rho}
c_r^{-1} \rho(z)\leq  \rho (\zeta)\leq c_r\,  \rho (z), \qquad \zeta\in D\bigl(z,r\rho(z)\bigr).
\end{equation}

There are at least two related distances in $\C$ induced by $\mathcal F_\phi$. The first is given by the reproducing kernel $K_\phi$ (see \cite{AMc}*{p.128}):
\begin{equation}\label{eq:dk}
 d_K(z,\zeta)=\sqrt{1-\frac{|K_\phi(z,\zeta)|^2}{K_\phi(z,z) K_\phi(\zeta,\zeta)}},\qquad z,\zeta\in\C,
\end{equation}
while the second is the distance $d_B$ associated to the Bergman metric \begin{equation}\label{def: Bergman metric}
\Delta \log K_\phi (z,z) dz \otimes d\bar z.
\end{equation}
According to Proposition 8 (or 9) in \cite{ARSW}, the Bergman distance is the so-called \textit{inner distance} of $d_K$, that is, 
\[
 d_B(z,\zeta)=\inf\bigl\{\ell_{d_k}(\gamma)\, :\, \text{$\gamma$ curve joining $z$ to $\zeta$}\bigr\} ,
\]
where, if $\gamma: [0,1]\longrightarrow \C$, 
\[
 \ell_{d_k}(\gamma)=\sup\left\{\sum_{i=1}^{n} d_K\bigl(\gamma(t_{i-1}), \gamma(t_{i})\bigr)\, :\, 0=t_0<\cdots<t_n=1\right\}.
\]
In particular,
 \begin{equation}\label{eq:d_k-less-d_B}
  d_K(z,\zeta)\leq d_B(z,\zeta)\, , \qquad z,\zeta\in\C.
 \end{equation}

Next, we recall some well-known estimates for the reproducing kernel (see \cite{BMOC}*{Proposition 14} and the references therein).

\begin{proposition}\label{prop:Jerry}
 There exist positive constants $C$, $\varepsilon$ and $\eta$, depending only on $\phi$, such that,
 \begin{itemize}
  \item [(a)] $|K_\phi(z,\zeta)|\leq C\, e^{\phi(z)+\phi(\zeta)} e^{-d_B^\varepsilon(z,\zeta)}$  for all $z,\zeta\in\C$,
  \item [(b)] $|K_\phi(z,\zeta)|\geq C\, e^{\phi(z)+\phi(\zeta)}$ for all $\zeta\in D(z,\eta\rho(z))$,
  \item [(c)] $C^{-1}/\rho^2(z)\leq \Delta\log K_\phi (z,z)\leq C/\rho^2(z)$.
 \end{itemize}
\end{proposition}

In \cite{BMOC},  estimate (a) is given in terms of the distance $d_\phi$ associated to the metric
$\rho^{-2}(z) dz\otimes d\bar z$ which, by (c), is equivalent to $d_B$.

Another property to be used in the proof of Theorem~\ref{main} is the following.

\begin{lemma}\label{lemma:Jerry}  (\cite{BMOC}*{Lemma 12}; \cite{MOC}*{Lemma 2.7})
 For any $\varepsilon>0$ and  $k\geq 0$ 
 \[
  \sup_{z\in\C}\int_{\C} \Bigl(\frac{|\zeta-z|}{\rho(z)}\Bigr)^k e^{- d_B^\varepsilon(z,\zeta)}\, \frac{dm(\zeta)}{\rho^2(\zeta)}<+\infty.
 \]
\end{lemma}

It is not difficult to see (\cite{MMOC}*{Lemma 4}) that for every $r>0$ there exists $C>0$ such that 
\begin{equation}\label{eq:dphi}
 C^{-1} \frac{|z-\zeta|}{\rho(z)}\leq d_B(z,\zeta)\leq C\, \frac{|z-\zeta|}{\rho(z)},\qquad \zeta\in D(z,r\rho (z)).
\end{equation}

The distance $d_K$ satisfies similar estimates, which will be important in the proof of Theorem~\ref{main}.

\begin{lemma}\label{lem:d_k}
(a) For every $r>0$ there exists $C>0$ depending only on $\phi$ such that 
\[
 d_K(z,\zeta)\leq C\, \frac{|z-\zeta|}{\rho (z)},\qquad \zeta\in D(z,r\rho (z)).
\]

(b) There exist $\eta>0$ and a positive constant $C$ such that
\[
 d_K(z,\zeta)\geq C\, \frac{|z-\zeta|}{\rho (z)},\qquad \zeta\in D(z,\eta\rho (z)).
\]
\end{lemma}

\begin{proof} (a) Apply \eqref{eq:d_k-less-d_B} and the upper estimate in \eqref{eq:dphi}.
 
(b)  By \cite{ARSW}*{p. 30}, the distance $d_K(z,\zeta)$ is the solution to the extremal problem
\begin{equation*}\label{eq:extremal}
 d_K(z,\zeta)=\frac{\sup\{|F(\zeta)| : F\in\mathcal F_\phi, \ \|F\|_\phi\leq 1,\ F(z)=0\}}{\sqrt{K_\phi(\zeta,\zeta)}}.
\end{equation*}
Take, fixed $z\in\C$, the function
\[
 F_z(\zeta)=\frac{\zeta-z}{\rho(z)}\frac{K_\phi(\zeta,z)}{\sqrt{K_\phi(z,z)}}.
\]
By Proposition~\ref{prop:Jerry}(a) and \eqref{eq:kernel-diagonal},
\begin{align*}
 \|F_z\|_\phi^2&\simeq \int_{\C}\frac{|\zeta-z|^2}{\rho^2(z)}\frac {|K_\phi(\zeta,z)|^2}{K_\phi(z,z)}\, e^{-2\phi(\zeta)}\frac{dm(\zeta)}{\rho^2(\zeta)}
 \lesssim \int_{\C}\frac{|\zeta-z|^2}{\rho^2(z)} e^{-2 d_B^\varepsilon(z,\zeta)} \frac{dm(\zeta)}{\rho^2(\zeta)},
\end{align*}
hence,  by Lemma~\ref{lemma:Jerry}, $\|F_z\|_\phi\lesssim 1$.

Thus, by  \eqref{eq:kernel-diagonal} once more, 
\begin{align*}
 d_K(z,\zeta) &\succsim \frac{|F_z(\zeta)|}{\sqrt{K_\phi(\zeta,\zeta)}}=\frac{|\zeta-z|}{\rho(z)}
 \frac{|K_\phi(\zeta,z)|}{\sqrt{K_\phi(z,z)} \sqrt{K_\phi(\zeta,\zeta)}}\\
 &\simeq \frac{|\zeta-z|}{\rho(z)}\, |K_\phi(\zeta,z)| e^{-\phi(z)} e^{-\phi(\zeta)}.
\end{align*}
Then we conclude with Proposition~\ref{prop:Jerry} (b).
\end{proof}

\subsection{Proof of Theorem~\ref{main}. Main Proposition}\label{subsec:proof-main}
We start with an easy remark. If 
\[
\int_{\C}\frac{dm(z)}{\rho^\gamma(z)}<+\infty, 
\] 
for some $\gamma>0$, then 
\begin{equation}\label{eq:rho-big}
 \liminf_{z\to+\infty} \rho(z)=+\infty.
\end{equation}

In fact, assume that \eqref{eq:rho-big} does not hold and take $C>0$ and $\{z_k\}_k\to\infty$ with $\rho(z_k)\leq C$. Assume furthermore that the disks $D(z_k,1)$ are pairwise disjoint. By the Lipschitz property \eqref{eq:rho-lip} we see that, for $k\geq 1$,
\[
 \rho(z)\leq C+1,\qquad z\in D(z_k,1).
\]
This leads to contradiction:
\[
 \int_{\C}\frac{dm(z)}{\rho^\gamma(z)}\geq\sum_k \int_{D(z_{k},1)}\frac{dm(z)}{(C+1)^\gamma}=\infty.
\]

From this point onward, we assume that \eqref{eq:rho-big} holds.
In view of \eqref{eq:1st-intensity}, this corresponds to the cases where the first intensity is asymptotically smaller than the area measure. In such situation, separation is more likely, and deciding whether $\Lambda_\phi$ is almost surely separated is more laborious. As a matter of fact, it is easy to see that if \eqref{eq:rho-big} does not hold, then $\Lambda_\phi$ is a.s. not separated (see Lemma~\ref{rem:rho-no-infinity} at the end of Subsection~\ref{subsec:conclusion}).

We begin the proof by considering a standard partition of the plane and studying the probability that two or more points fall in the same cell. 

Let
\[
T_{n,k}:= \Big\{z=re^{i\theta} \in \C: n-1\leq r<n, \, \frac{\theta}{2\pi} \in \big[\frac{k-1}{n}, \frac{ k}{n}\big)\Big\}, \qquad n\geq 1,\ k=1,\dots, n
\]
and let $z_{n,k}=(n-1/2)\, e^{i 2\pi\frac{k-1/2}n}$ denote its ``center''. 
Observe that, by \eqref{eq:rho-big}, for any $\eta>0$ there exists $n_0$ such that for $n\geq n_0$ and $k$, the cell $T_{n,k}$ is included in the disk $D\bigl(z_{n,k}, \eta\rho(z_{n,k})\bigr)$, so that $\rho(z)\simeq \rho(z_{n,k})$ for $z\in T_{n,k}$.

Consider the corresponding counting random variables
\[
 X_{n,k}=\#(\Lambda_{\phi}\cap T_{n,k}).
\]
Observe that, by \eqref{eq:EB} and \eqref{eq:rho},

\begin{equation}\label{expectation-xnk}
 \mathbb E[X_{n,k}]\simeq \int_{T_{n,k}} \frac{dm(z)}{\rho^2(z)}\simeq \frac 1{\rho^2(z_{n,k})}.
\end{equation}

A specific feature of any determinantal process $\Lambda$ is that the counting random variable $N(B)=\#(\Lambda\cap B)$ can be expressed as a sum of independent Bernoulli variables $\xi_j$. Moreover, the parameters $\lambda_j$ of the random variables $\xi_j$ are precisely the eigenvalues of the  restriction operator on $B$ (see e.g. \cite{HKPV}*{Theorem 4.5.3}). Hence, taking $\Lambda=\Lambda_\phi$ and $B=T_{n,k}$, we have
\begin{equation*}\label{eq:xnk-sum-bernoullis}
X_{n,k} =\sum_{j=1}^\infty \xi_j,
\end{equation*}
where $\xi_j\sim\text{Bernoulli}(\lambda_j)$ and $\lambda_j$ is the $j^{\text{th}}$ eigenvalue (arranged in decreasing order) of the restriction operator
\begin{equation}\label{eq:restriction-operator}
 Tf(z)=\int_{T_{n,k}} f(\zeta) K_{\phi}(z,\zeta) d\mu_{\phi}(\zeta),\qquad f\in\mathcal F_{\phi}.
\end{equation}

The key result in the proof of Theorem~\ref{main} is the estimate of the probability that two or more points fall in the same cell $T_{n,k}$, which, as we shall see soon, is comparable to  $\rho^{-6}(z_{n,k})$. 
We begin with the following general representation of this probability.

\begin{lemma}\label{lemma:crucial}
 Let $T_{n,k}$ and $X_{n,k}$ as above. Then,
 \begin{align*}
  \P(X_{n,k}\geq 2)=\int_{T_{n,k}}\int_{T_{n,k}} [K_\phi(z,z) K_\phi(\zeta,\zeta)- |K_\phi(z,\zeta)|^2]\, d\mu_\phi(z)\,d\mu_\phi(\zeta)
  + O\bigl((\mathbb E[X_{n,k}])^3\bigr).
 \end{align*}
\end{lemma}

\begin{remark}\label{rem:EXX-1}
(a) As we shall see in the proof, this identity is valid in great generality. Given a determinantal point process $\Lambda$ with correlation kernel $K(z,\zeta)$ and reference measure $\mu$, and given $B\subset\C$ and its counting variable $N(B)=\#(\Lambda\cap B)$, the analogue for $\P(N(B)\geq 2)$ also holds. Of course, this is only of interest when $(\mathbb E[N(B)])^3$ is much smaller than the double integral.

(b) The identity in the statement can also be rewritten in terms of the distance $d_K$ defined in \eqref{eq:dk} and the first intensity of the process:
\[
 \P(X_{n,k}\geq 2)=\int_{T_{n,k}}\int_{T_{n,k}} d_K^2(z,\zeta) K_\phi(z,z)\, d\mu_\phi(z)\ K_\phi(\zeta,\zeta)\, d\mu_\phi(\zeta)
  + O\bigl((\mathbb E[X_{n,k}])^3\bigr).
\]

(c) The double integral above is exactly $\mathbb E[X_{n,k}(X_{n,k}-1)]$ (see \cite{HKPV}*{(1.2.4)}). Thus, the result can be rephrased as
\[
 \P(X_{n,k}\geq 2)=\mathbb E[X_{n,k}^2]-\mathbb E[X_{n,k}]+ O\bigl((\mathbb E[X_{n,k}])^3\bigr).
\]
\end{remark}

\begin{proof}
 By definition,
 \begin{align*}
   \P(X_{n,k}\geq 2)&=1- \P(X_{n,k}= 0)- \P(X_{n,k}= 1)=1-\prod_{j=1}^\infty(1-\lambda_j)-\sum_{l=1}^\infty \lambda_l \prod_{j\neq l}(1-\lambda_j)\\
   &=1-\exp\Bigl(-\sum_{j=1}^\infty\log\bigl(\frac 1{1-\lambda_j}\bigr)\Bigr)-
   \exp\Bigl(-\sum_{j=1}^\infty\log\bigl(\frac 1{1-\lambda_j}\bigr)\Bigr)\Bigl(\sum_{l=1}^\infty\frac{\lambda_l}{1-\lambda_l}\Bigr).
 \end{align*}
Observe that all $\lambda_j$ are very close to 0, because of condition \eqref{eq:rho-big}, and estimate \eqref{expectation-xnk}:
\[
 \E[X_{n,k}]=\text{trace}(T)=\sum_{j=1}^\infty\lambda_j\simeq \frac 1{\rho^2(z_{n,k})}.
\]

Plugging in the power series 
\begin{align*}
 & e^{-x}=\sum_{n=0}^\infty(-1)^{n+1}\frac{x^n}{n!}=1-x+\frac{x^2}2+\cdots \\
& \log\bigl(\frac 1{1-x}\bigr)=\sum_{n=1}^\infty\frac{x^n}n=x+\frac{x^2}{2}+\cdots \\
& \frac x{1-x}=\sum_{n=1}^\infty x^n=x+x^2+\cdots
\end{align*}
and keeping the terms up to order 2 in the expression above, we get
\begin{align*}
 \P(X_{n,k}\geq 2)&= \sum_{j=1}^\infty\log\bigl(\frac 1{1-\lambda_j}\bigr)-\frac 12\Bigl(\sum_{j=1}^\infty\log\bigl(\frac 1{1-\lambda_j}\bigr)\Bigr)^2- \\
 &\qquad -\Bigl(\sum_{l=1}^\infty\frac{\lambda_l}{1-\lambda_l}\Bigr)\left(1-\sum_{j=1}^\infty\log\bigl(\frac 1{1-\lambda_j}\bigr)\right)+O\Bigl(\bigl(\sum_{j=1}^\infty\lambda_j\bigr)^3\Bigr)\\
 &= \sum_{j=1}^\infty (\lambda_j+\frac{\lambda_j^2}2)-\frac 12 \Bigl(\sum_{j=1}^\infty (\lambda_j+\frac{\lambda_j^2}2)\Bigr)^2
 -\Bigl(\sum_{l=1}^\infty (\lambda_l+\lambda_l^2)\Bigr)-\\
 &\qquad\qquad -\Bigl(\sum_{l=1}^\infty(\lambda_l+\lambda_l^2)\Bigr)\Bigl(1-\sum_{j=1}^\infty (\lambda_j+\frac{\lambda_j^2}2)\Bigr)+O\Bigl(\bigl(\sum_{j=1}^\infty\lambda_j\bigr)^3\Bigr)\\
 &=\sum_{j,l=1}^\infty \lambda_l\lambda_j-\sum_{l=1}^\infty\lambda_l^2+O\Bigl(\bigl(\sum_{j=1}^\infty\lambda_j\bigr)^3\Bigr)\\
& =\bigl(\sum_{l=1}^\infty\lambda_l\bigr)^2-\sum_{l=1}^\infty\lambda_l^2+O\Bigl(\bigl(\sum_{j=1}^\infty\lambda_j\bigr)^3\Bigr).
\end{align*}

By the trace formula, 
\[
 \text{trace}(T)=\sum_{j=1}^\infty \lambda_j=\int_{T_{n,k}} K_{\phi}(\zeta,\zeta) d\mu_{\phi}(\zeta),
\]
and by the spectral theorem for integral Hilbert-Schmidt operators,
\[
 \sum_{j=1}^\infty \lambda_j^2=\int_{T_{n,k}} \int_{T_{n,k}} |K_{\phi}(z, \zeta)|^2 d\mu_{\phi}(z) d\mu_{\phi}(\zeta).
\]
Thus
\begin{align*}
\bigl(\sum_{j=1}^\infty \lambda_j\bigr)^2-\sum_{j=1}^\infty \lambda_j^2
=\int_{T_{n,k}}\int_{T_{n,k}} \left[K_{\phi}(z,z) K_{\phi}(\zeta,\zeta)-|K_{\phi}(z,\zeta)|^2\right] d\mu_{\phi}(z) d\mu_{\phi}(\zeta) ,
\end{align*}
as desired.
\end{proof}

An important consequence of Lemma~\ref{lemma:crucial} is the following.

\begin{proposition}\label{prop:traces}
 There exists $n_0=n_0(\phi)\in\N$ such that for all $n\geq n_0$ and $k=1,\dots, n$
 \[
  \P(X_{n,k}\geq 2)\simeq \frac {1}{\rho^6(z_{n,k})},
 \]
 where the implicit constants do not depend on $n$ and $k$.
\end{proposition}

\begin{proof}
We start with the formula of Remark~\ref{rem:EXX-1}(b). 
Observe that, by condition \eqref{eq:rho-big}, for any $\eta>0$ there exists $n_0$ such that $T_{n,k}\subset D(z,\eta\rho(z))$ for all $z\in T_{n,k}$, $n\geq n_0$ and $k=1,\dots,n$. Then, by Lemma~\ref{lem:d_k}, 
\[
 d_K(z,\zeta)\simeq \frac{|z-\zeta|}{\rho(z)}, \qquad z,\zeta\in T_{n,k}.
\]
Therefore
\begin{align*}
 \P(X_{n,k}\geq 2)&
  \simeq\int_{T_{n,k}}\int_{T_{n,k}}\frac{|z-\zeta|^2}{\rho^2(z)}\,  K_\phi(z,z)\, d\mu_\phi(z)\ K_\phi(\zeta,\zeta)\, d\mu_\phi(\zeta)
  + O\Bigl((\mathbb E[X_{n,k}])^3\Bigr).
\end{align*}
By \eqref{eq:1st-intensity}, \eqref{expectation-xnk} and \eqref{eq:rho},
we get finally
\begin{align*}
 \P(X_{n,k}\geq 2)&\simeq \int_{T_{n,k}}\int_{T_{n,k}}\frac{|z-\zeta|^2}{\rho^2(z)}\, \frac{dm(z)}{\rho^2(z)}\frac{dm(\zeta)}{\rho^2(\zeta)}+ O\Bigl(\frac {1}{\rho^6 (z_{n,k})}\Bigr)\\
&\simeq\frac {1}{\rho^6 (z_{n,k})}\left(\int_{T_{n,k}}\int_{T_{n,k}}|z-\zeta|^2\, dm(z)\, dm(\zeta) +1\right)\simeq \frac {1}{\rho^6 (z_{n,k})}.
\end{align*}
\end{proof}

\begin{remark}\label{rem:lambda12} 
The lower estimate in Proposition~\ref{prop:traces} can be proved in a more direct way. Since the expected value of $X_{n,k}=\sum_j \xi_j$ is close to zero, the probability that two or more points fall in $T_{n,k}$ is concentrated in the events $\{\xi_1=1\}$ and $\{\xi_2=1\}$ (since $\xi_1,\xi_2$ are the two Bernoulli variables with larger parameter). In any case, by the independence of $\xi_1$ and $\xi_2$, we have the estimate
\[
\P(X_{n,k}\geq 2)\geq \P\bigl((\xi_1=1)\cap (\xi_2=1) \bigr)=\lambda_1\, \lambda_2.
\]
By definition
\[
 \lambda_1=\sup_{\substack{f\in \mathcal F_{\phi}\\ \|f\|_{\phi}=1}} \langle Tf,f \rangle_{\phi},
\]
and, by \eqref{eq:restriction-operator} and the reproducing property of $K_{\phi}(z,\zeta)$,
\[
 \langle Tf,f \rangle_{\phi}=\int_{T_{n,k}} |f(z)|^2 d\mu_{\phi}(z).
\]

Taking as $f$ the normalized reproducing kernel for the center of $T_{n,k}$,
\[
 \mathfrak K_{n,k}(z)=\frac{K_{\phi}(z,z_{n,k})}{\sqrt{K_{\phi}(z_{n,k}, z_{n,k})}}
\]
we see that 
\[
 \lambda_1\geq \int_{T_{n,k}} \frac{|K_{\phi}(z,z_{n,k})|^2}{K_{\phi}(z_{n,k}, z_{n,k})} \, d\mu_{\phi}(z).
\]
By Proposition~\ref{prop:Jerry} (b), 
\begin{equation}\label{eq:lambda1}
 \lambda_1\succsim \int_{T_{n,k}} K_\phi(z,z)\, d\mu_\phi(z)\simeq \int_{T_{n,k}} \frac{dm(z)}{\rho^2(z)}\simeq\frac 1{\rho^2(z_{n,k})}.
\end{equation}

In order to estimate $\lambda_2$ we use the formula
\[
 \lambda_2=\sup_{\substack{\text{dim}(F)=2 \\ F\subset \mathcal F_{\phi}}}\ \inf_{g\in F} \frac {\langle Tg,g \rangle_{\phi}}{\|g\|_{\phi}^2}.
\]
Take $F\subset \mathcal F_{\phi}$ to be the subspace generated by 
\[
 g_1(z)=\mathfrak K_{n,k}(z)\qquad\text{and}\qquad g_2(z)=\frac {z-z_{n,k}}{\rho(z_{n,k})}\, \mathfrak K_{n,k}(z).
\]
Proposition~\ref{prop:Jerry} (a), together with \eqref{eq:kernel-diagonal} yield, for some $\varepsilon>0$,
\[
 \|g_2\|_{\phi}^2=\int_{\C} \frac {|z-z_{n,k}|^2}{\rho^2(z_{n,k})}\, \frac{|K_{\phi}(z,z_{n,k})|^2}{K_{\phi}(z_{n,k}, z_{n,k})} \, d\mu_{\phi}(z)\lesssim \int_{\C} \frac {|z-z_{n,k}|^2}{\rho^2(z_{n,k})}\,e^{-d_B^\varepsilon(z,z_{n,k})} \frac{dm(z)}{\rho^2(z)}.
\]
By Lemma~\ref{lemma:Jerry} with $k=2$, we see that $\|g_2\|_{\phi}\lesssim 1$.

Also, since $g_2$ vanishes on $z_{n,k}$, the functions $g_1$ and $g_2$ are orthogonal in $\mathcal F_{\phi}$. Then, for any function
$g=a g_1+b g_2\in F$, $ a,b\in\C$, 
one has 
\[
\|g\|_{\phi}^2=|a|^2 \|g_1\|_{\phi}^2+|b|^2\|g_2\|_{\phi}^2\simeq |a|^2+|b|^2.
\]

On the other hand, again by the estimates of Proposition~\ref{prop:Jerry} (a) and (b),
\begin{align*}
 \langle Tg,g \rangle_{\phi}&\simeq \int_{T_{n,k}} \Bigl|a+b \frac {\zeta-z_{n,k}}{\rho(z_{n,k})}\Bigr|^2 
 \frac{|K_{\phi}(\zeta,z_{n,k})|^2}{K_{\phi}(z_{n,k},z_{n,k})}\ d\mu_{\phi}(\zeta)\\
 &\simeq  \int_{T_{n,k}} \bigl|a+b \frac {\zeta-z_{n,k}}{\rho(z_{n,k})}\bigr|^2 \ 
  \frac{dm(\zeta)}{\rho^2(\zeta)}\\
& \simeq \frac 1{\rho^2(z_{n,k})} \int_{T_{n,k}} \bigl|a+b \frac {\zeta-z_{n,k}}{\rho(z_{n,k})}\bigr|^2\ dm(\zeta).
\end{align*}
Let $D_{n,k}=D(z_{n,k},1/4)$ and notice that $D_{n,k}\subset T_{n,k}$. Since $1$ and $\zeta-z_{n,k}$ are orthogonal in $L^2(D_{n,k}, dm)$, we see that
\begin{align*}
  \langle Tg,g \rangle_{\phi}& \succsim  \frac 1{\rho^2(z_{n,k})} \int_{D_{n,k}} \bigl|a+b \frac {\zeta-z_{n,k}}{\rho(z_{n,k})}\bigr|^2\ dm(\zeta)\\
  &\simeq \frac 1{\rho^2(z_{n,k})} \Bigl(|a|^2 +\frac{|b|^2}{\rho^2(z_{n,k})} \int_{D_{n,k}} |\zeta-z_{n,k}|^2 dm(\zeta)\Bigr)\\
  &\simeq \frac 1{\rho^2(z_{n,k})} \bigl(|a|^2+\frac{|b|^2}{\rho^2(z_{n,k})}\bigr).
\end{align*}
With this
\[
 \lambda_2\geq \inf_{g\in F} \frac {\langle Tg,g \rangle_{\phi}}{\|g\|_{\phi}^2} \succsim\inf_{a,b\in\C}\frac {\frac 1{\rho^2(z_{n,k})} \bigl(|a|^2+\frac{|b|^2}{\rho^2(z_{n,k})}\bigr)}{|a|^2+|b|^2}=\frac 1{\rho^4(z_{n,k})},
\]
which together with \eqref{eq:lambda1} yields finally
\[
 P(X_{n,k}\geq 2)\geq\lambda_1 \lambda_2\succsim \frac 1{\rho^2(z_{n,k})} \frac 1{\rho^4(z_{n,k})}=\frac {1}{\rho^6(z_{n,k})}.
\]
\end{remark}

\subsection{Proof of Theorem~\ref{main}. Conclusion}\label{subsec:conclusion}
A direct consequence of Proposition~\ref{prop:traces} is the following.

\begin{proposition}\label{prop:main-partition}
Let $\phi$ satisfy \eqref{eq:rho-big} and let $X_{n,k}$ be as above. Then
\[
 \P\bigl(X_{n,k}\geq 2\ \text{infinitely often}\bigr)=
 \begin{cases}
  0\quad  \text{if $\displaystyle{\int_{\C}}\dfrac{dm(z)}{\rho^6(z)}<+\infty$}\\
  1 \quad  \text{if $\displaystyle{\int_{\C}}\dfrac{dm(z)}{\rho^6(z)}=+\infty$.}
 \end{cases}
\]
\end{proposition}

\begin{proof}
 Assume first that $\int_{\C}\frac{dm(z)}{\rho^6(z)}<+\infty$. By Proposition~\ref{prop:traces}
 \[
  \sum_{n,k} \P(X_{n,k}\geq 2)\simeq \sum_{n,k} \frac {1}{\rho^6(z_{n,k})} \simeq \sum_{n,k}\int_{T_{n,k}}\frac{dm(z)}{\rho^6(z)}= \int_{\C}\frac{dm(z)}{\rho^6(z)} <\infty,
 \]
 and the result follows from the first Borel-Cantelli lemma \cite{Billingsley95}*{Theorem 4.3}.
 
 Assume now that $\int_{\C}\frac{dm(z)}{\rho^6(z)}=+\infty$. Observe that if 
\[
 E=\{ X_{n,k}\geq 2\quad\text{infinitely often}\},
\]
then 
\begin{align*}
 E^c&=\{\exists n_0\in\N : \forall n\geq n_0,\ \forall k=1,\dots, n,\quad X_{n,k}\leq 1\}
 =\bigcup_{n_0=1}^\infty \bigcap_{n\geq n_0}\bigcap_{k=1}^n (X_{n,k}\leq 1).
\end{align*}
Thus, in order to see that $\P(E^c)=0$ it will be enough to see that for all $n_0\in\N$ 
\begin{equation}\label{eq:event-n0}
 \P\bigl(\bigcap_{n\geq n_0}\bigcap_{k=1}^n (X_{n,k}\leq 1)\bigr)=0.
\end{equation}
 According to a general estimate by S. Ghosh \cite{Gh}*{Theorem 1.4}, for every collection of disjoint Borel sets $\{B_j\}_{j=1}^N$ and integers $\{m_j\}_{j=1}^N$
 \begin{equation}\label{eq:Ghosh}
  \P\bigl(\cap_{j=1}^N (N_{\phi}(B_j)\leq m_j)\bigr)\leq \prod_{j=1}^N \P\bigl( (N_{\phi}(B_j)\leq m_j\bigr),
 \end{equation}
 where $N_{\phi}(B)=\#(\Lambda_{\phi}\cap B )$ indicates the counting function on $B$.
 
Taking any $N>n_0$ and applying this estimate to all $T_{n,k}$ with $ n_0\leq n\leq N$, we see that
\[
  \P\bigl(\bigcap_{n=n_0}^N \bigcap_{k=1}^n (X_{n,k}\leq 1)\bigr)
  \leq \prod_{n=n_0}^N \prod_{k=1}^n \P\bigl(X_{n,k}\leq 1\bigr).
\]
Letting $N\to\infty$ we get
\[
 \P\bigl(\bigcap_{n\geq n_0} \bigcap_{k=1}^n (X_{n,k}\leq 1)\bigr)\leq \prod_{n=n_0}^\infty \prod_{k=1}^n \P\bigl(X_{n,k}\leq 1\bigr),
\]
hence we shall have \eqref{eq:event-n0} as soon as we prove that for all $n_0\in\N$
\begin{equation}\label{eq:product}
 \prod_{n=n_0}^\infty \prod_{k=1}^n \P\bigl(X_{n,k}\leq 1\bigr)=0.
\end{equation}

Denoting $\varepsilon_{n,k}=\P(X_{n,k}\geq 2)$ and recalling that, by Proposition~\ref{prop:traces},
$\varepsilon_{n,k}\simeq \rho^{-6}(z_{n,k})$, we have
\begin{align*}
  \prod_{n=n_0}^\infty \prod_{k=1}^n \P\bigl(X_{n,k}\leq 1\bigr)&=\exp\Bigl(-\sum_{n=n_0}^\infty \sum_{k=1}^n \log\bigl(\frac 1{1-\varepsilon_{n,k}}\bigr)\Bigr)\\
 &\simeq \exp\Bigl(-\sum_{n=n_0}^\infty\sum_{k=1}^n  \varepsilon_{n,k}\Bigr)\simeq \exp\Bigl(- \sum_{n=n_0}^\infty \sum_{k=1}^n
 \frac 1{\rho^6(z_{n,k})}\Bigr).
\end{align*}
Since the series 
\[
 \sum_{n=1}^\infty \sum_{k=1}^n \frac 1{\rho^6(z_{n,k})}\simeq \int_{\C}\frac{dm(z)}{\rho^6(z)}
\]
diverges, we have \eqref{eq:product}, as desired.
\end{proof}

Once Proposition~\ref{prop:main-partition} is settled, we finish the proof of Theorem~\ref{main} for the case \eqref{eq:rho-big} with a standard argument dating back, at least, to the proof of Theorem 2 in \cite{Coch}.

Assume first that $\int_{\C}\frac{dm(z)}{\rho^6(z)}<\infty$.  
All the arguments above can be applied analogously to the shifted regions ($ n\geq 1,\ k=1,\dots, n$):
\begin{align*}
 &\widetilde T_0=\{z\in\C : |z|<1/2\}=D(0,1/2),\\
 &\widetilde T_{n,k}=\Big\{z=re^{i\theta} \in\C : n-\frac 12\leq r < n+\frac 12,\, \frac{\theta}{2\pi} \in \big[\frac{k-1/2}{n}, \frac{ k+1/2}{n}\big)\Big\}, 
\end{align*}
and the corresponding random variables $\widetilde X_{n,k}=\#(\Lambda_{\phi} \cap \widetilde T_{n,k})$. In particular, Proposition~\ref{prop:main-partition} holds equally if the random variables $X_{n,k}$ are replaced by the similar random variables $\tilde X_{n,k}$.

Consider then the events
\begin{align*}
 E=\{X_{n,k}\geq 2\ \text{infinitely often}\},\quad \widetilde E=\{\widetilde X_{n,k}\geq 2\ \text{infinitely often}\},
\end{align*}
for which $\mathbb P(E^c)=\mathbb P(\widetilde E^c)=1$.
Under the event $E^c\cap\widetilde E^c$, for all but finitely many $n$, $k$, we have at most one point in the regions $T_{n,k}$, $\widetilde T_{n,k}$. Thus, the points of $\Lambda_{\phi}$ contained in these regions are separated by a fixed constant ($1/4$ will do). Since the finitely many regions where $X_{n,k}$, $\widetilde X_{n,k}$ could be bigger than 1 contain at most a finite number of points of $\Lambda_{\phi}$, we deduce that $\Lambda_{\phi}$ is almost surely separated. 

Notice that we actually prove that $\Lambda_\phi$ is the union of a finite sequence and a sequence which is separated by a constant comparable to the size of the cells $T_{n,k}$.

For the case $\int_{\C}\frac{dm(z)}{\rho^6(z)}=\infty$ we can repeat the above arguments to a grid of size $1/l$, $l\geq 1$. Let
\[
 T_{n,k}^l:= \Big\{z=r e^{i\theta} \in \C: \frac{n-1}l\leq r<\frac nl, \, \frac{\theta}{2\pi} \in \big[\frac{k-1}{ln}, \frac{ k}{ln}\big)\Big\}, \qquad n\geq 1,\ k=1,\dots, ln.
\]
Then, letting $X_{n,k}^{(l)}=\#(\Lambda_{\phi}\cap T_{n,k}^l)$, we see, in the same way as in Proposition~\ref{prop:main-partition}, that the events
\[
 E^l=\{X_{n,k}^{(l)}\geq 2\ \text{infinitely often}\}
\]
have all probability 1. 
Under the event $E^l$ there are infinitely many couples $\lambda_k, \lambda_j\in\Lambda_{\phi}$ at a distance smaller than $1/l$. Therefore, under the event $\cap_{l=1}^\infty E^l$, which still has probability 1, the sequence $\Lambda_{\phi}$ is not separated.

It remains to settle the case where \eqref{eq:rho-big} is not satisfied. Similar arguments to those is Remark~\ref{rem:lambda12} show that the process is then a.s not separated.

\begin{lemma}\label{rem:rho-no-infinity}
	If condition \eqref{eq:rho-big} does not hold, then $\Lambda_\phi$ is almost surely not separated.
\end{lemma}

\begin{proof}
	 Assume that there exist $C>0$ and $\{z_k\}_{k=1}^\infty\subset\C$ such that $\{z_k\}_{k}\to \infty$ and $\rho(z_k)\leq C$ for all $k$. 
	 For $l\in\mathbb N$ consider the disks $D_k^l=D\bigl(z_k,\rho(z_k)/l\bigr)$ and their associated counting variables $X_k^{(l)}=\#(\Lambda_\phi\cap D_k^l)$. As seen in Remark~\ref{rem:lambda12},
	 \[
	  \mathbb P (X_k^{(l)}\geq 2)\geq \lambda_1\lambda_2,
	 \]
	 where $X_k^{(l)}=\sum_{j=1}^\infty\xi_j$ is a sum of independent variables $\xi_j\sim\text{Bernoulli}(\lambda_j)$ and $\{\lambda_j\}_{j=1}^\infty$ are the eigenvalues, arranged in decreasing order, of the restriction operator
	 \[
	 Tf(z)=\int_{D_k^l} f(\zeta) K_{\phi}(z,\zeta) d\mu_{\phi}(\zeta),\qquad f\in\mathcal F_{\phi}.
	 \]
	 As before, taking the normalized reproducing kernel $\mathfrak K_k$ at $z_k$, and using the estimates of Proposition~\ref{prop:Jerry}(b), we see that for $l$ big enough
	 \[
	  \lambda_1\geq \langle T \mathfrak K_k,\mathfrak K_k\rangle\simeq \int_{D_k^l}\frac{dm(z)}{\rho^2(z)}\simeq \frac 1{l^2}.
	 \]
	Considering the space $F$ in $\mathcal F_\phi$ generated by $g_1=\mathfrak K_k$ and $g_2(z)=\frac{z-z_k}{\rho(z_k)} \mathfrak K_k(z)$, we have now
	\begin{align*}
	 \lambda_2 &\geq \inf_{g\in F} \frac {\langle Tg,g \rangle_{\phi}}{\|g\|_{\phi}^2} 
	 \succsim\inf_{a,b\in\C}\frac {\int_{D_k^l}\bigl|a+b\frac{z-z_k}{\rho(z_k)}\bigr|^2 \frac{dm(z)}{\rho^2(z)}}{|a|^2+|b|^2}\\
	 & \simeq \frac 1{\rho^2(z_k)} \inf_{a,b\in\C}\frac{|a|^2 |D_k^l|+|b|^2\int_{D_k^l} \frac{|z-z_k|^2}{\rho^2(z_k)} dm(z)}{|a|^2+|b|^2}\\
	 &\simeq \frac 1{\rho^2(z_k)} \inf_{a,b\in\C}\frac{|a|^2 \rho^2(z_k)/l^2+|b|^2 \rho^2(z_k)/l^4}{|a|^2+|b|^2} \simeq \frac 1{l^4}.
	\end{align*}
	With this
	\[
	 \mathbb P (X_k^{(l)}\geq 2) \succsim \frac 1{l^6},
	\]
	and using again Ghosh's estimate \eqref{eq:Ghosh},
	\[
	 \mathbb P\bigl[\cap_{k=1}^\infty (X_k^{(l)}\leq 1)\bigr]\leq\prod_{k=1}^\infty \mathbb P(X_k^{(l)}\leq 1)=
	 \prod_{k=1}^\infty \bigl(1-\mathbb P(X_k^{(l)}\geq 2)\bigr)=0.
	\]
	Thus, the events $E_l=\{X_k^{(l)}\geq 2\ \text{infinitely often}\}$ have probability one, so that
	$E:=\cap_j E_l$
	 has also probability 1. Since $D_k^l\subset D(z_k,C/l)$ for all $k,l$, by the hypothesis on $\rho(z_k)$, we see that under the event $E$ the sequence $\Lambda_\phi$ is not separated.
\end{proof}

We conclude this subsection with a related result that of independent interest.

A sequence $\{\lambda_k\}_k\subset\C$ is separated with respect to a distance $d$ (or $d$-separated) if
\[
 \inf_{j\neq k} d(\lambda_j,\lambda_k)>0.
 \]
The arguments of Lemma \ref{rem:rho-no-infinity}, together with \eqref{eq:dphi}, show the following result.

\begin{proposition} \label{prop:d_b separation}
Let $d_B$ be the Bergman metric defined in Subection~\ref{subsec:prelim}. Then for any doubling subharmonic weight $\phi$
 \[
  \mathbb P\bigl(\Lambda_\phi\ \text{is $d_B$-separated}\bigr)=0.
 \]
 \end{proposition}
 
Since Lemma~\ref{lem:d_k} states that $d_K(z,\zeta) \simeq d_B(z,\zeta)$ whenever $\zeta \in D(z,\rho(z))$, the above results also holds for the distance $d_K$. 

\subsection{Proof of Corollary~\ref{cor:interpolation}}
Assume that $\int_{\C}\frac{dm(z)}{\rho^6(z)}<+\infty$. In particular, as observed at the beginning of Subsection~\ref{subsec:proof-main}, condition \eqref{eq:rho-big} is satisfied. 

Since $\Lambda_\phi$ is a.s. separated, it will be enough to see that the upper density of $\Lambda_\phi$ is zero almost surely. 
For any $N\in\N$ consider now the cells of size $N$,
\[
 T_{n,k}^N:= \Big\{z=re^{i\theta} \in \C: (n-1)\leq \frac rN < n, \ \frac{\theta}{2\pi N} \in \big[(k-1), k\big)\Big\}, \quad n\geq 1,\ k=1,\dots, n
\]
and their shifted versions
\begin{align*}
&\widetilde T_0^N=D(0,N/2)\\
&\widetilde T_{n,k}^N:= \Big\{z=re^{i\theta} \in \C: n-\frac 12\leq \frac rN < n+\frac 12, \ \frac{\theta}{2\pi N} \in \big[k-\frac 12, k+\frac 12\big)\Big\}.
\end{align*}
Consider also the associated counting variables 
\[
X_{n,k}^{(N)}=\#(\Lambda_\phi\cap T_{n,k}^N)\, ,\qquad  \widetilde X_{n,k}^{(N)}=\#(\Lambda_\phi\cap \widetilde T_{n,k}^N). 
\]
The same arguments used in the proof of Theorem~\ref{main} show that the events
\[
 E=\{X_{n,k}^{(N)}\geq 2\ \text{infinitely often}\},\quad \tilde E=\{\widetilde X_{n,k}^{(N)}\geq 2\ \text{infinitely often}\}
\]
have both probability $0$. Hence, after removing a finite number of points to $\Lambda_\phi$ we obtain a sequence in which the points are separated by a distance comparable to the size of the cells. A simple area estimate implies then that
\[
 D_+(\Lambda_\phi)=\limsup_{r\to\infty}\sup_{z\in\C}\frac{\#(\Lambda_\phi\cap D(z,r))}{\pi r^2}\lesssim \frac 1{N^2}.
\]
Since this is valid for any $N$, we deduce that almost surely $\Lambda_\phi$ is separated and $ D_+(\Lambda_\phi)=0$.

\section{Poisson process: proof of Theorem~\ref{thm:Poisson}}\label{sec:Poisson}
This is a straightforward adaptation of \cite{HaMa}*{Theorem 1.2}, which we briefly recall.

It will be enough to prove the analogue of Proposition~\ref{prop:main-partition} for the Poisson point process: if 
\[
 X_{n,k}^P=\#(\Lambda_{\phi}^P\cap T_{n,k}),
\]
then
 \begin{equation}\label{eq-sep}
  \mathbb P\bigl(X_{n,k}^P\geq 2\ \text{infinitely often}\bigr)=
  \begin{cases}
     1\ \text{if $\displaystyle{\int_{\C}}\dfrac{dm(z)}{\rho^4(z)}<+\infty$}\\
  0\ \text{if $\displaystyle{\int_{\C}}\dfrac{dm(z)}{\rho^4(z)}=+\infty$.}
  \end{cases}
 \end{equation}
 Once this is proved, the same arguments applied in the determinantal case to shifted and to thinner partitions of the plane 
 yield the result.

By definition of the Poisson process $\Lambda_{\phi}^P$, the parameter $\mathfrak p_{n,k}$ of the Poisson variable $X_{n,k}$ is
\begin{align*}
\mathfrak p_{n,k}= \mathbb E(X_{n,k}^P)&= \Var(X_{n,k}^P)= \int_{T_{n,k}} \frac{dm(z)}{\rho^2(z)}\simeq 
\frac 1{\rho^2(z_{n,k})}.
\end{align*}

Thus, 
\begin{align*}
  \mathbb P\bigl(X_{n,k}^P\geq 2\bigr)&=1-\mathbb P\bigl(X_{n,k}^P=0\bigr)-\mathbb P\bigl(X_{n,k}^P=1\bigr)=1-e^{-\mathfrak p_{n,k}}-e^{-\mathfrak p_{n,k}}\mathfrak p_{n,k}\\
  &=\frac{\mathfrak p_{n,k}^2}2\bigl(1+o(1)\bigr)\simeq \frac 1{\rho^4(z_{n,k})}.
\end{align*}
Therefore
\[
 \sum_{n=1}^\infty\sum_{k=1}^n \mathbb P\bigl(X_{n,k}^P\geq 2\bigr)\simeq \sum_{n=1}^\infty \sum_{k=1}^n \frac 1{\rho^4(z_{n,k})}\simeq 
\int_{\C}\dfrac{dm(z)}{\rho^4(z)},
\]
and \eqref{eq-sep} follows from the second Borel-Cantelli lemma (\cite{Billingsley95}*{Theorem 4.4}) since, by definition of Poisson process, the diferent variables $X_{n,k}^P$ are independent.


\bibliographystyle{plain}

\end{document}